\RequirePackage{fix-cm}
\documentclass[smallextended]{svjour3}       
\smartqed  
\usepackage{bm} 
\usepackage{graphicx,color}
\usepackage{amssymb}
\usepackage{amsmath}
\usepackage{ascmac}

\usepackage{algorithmicx}
\usepackage[ruled]{algorithm}
\usepackage{algpseudocode}

\usepackage[utf8]{inputenc}

\newcommand{\q}{\quad}

\usepackage{url}
\newcommand{\fdel}[1]{}
\newcommand{\finfo}[1]{}

\newcommand{\ffinfo}[1]{}

\begin{document}

\title{A Parallelizable Energy-Preserving Integrator MB4 and Its Application to Quantum-Mechanical Wavepacket Dynamics
\thanks{The present research was partially supported by JSPS KAKENHI Grant Numbers JP16KT0016, JP16K17550, JP17H02828, JP17K19966 and JP19KK02555. The present research was partially supported by Priority Issue 7 of the post-K project.}
}

\titlerunning{Application of the MB4 Energy-Preserving Integrator}        

\author{Tsubasa Sakai \and
Shuhei Kudo \and
Hiroto Imachi \and
Yuto Miyatake \and
Takeo Hoshi \and
Yusaku Yamamoto}


\institute{T.~Sakai and Y.~Yamamoto \at
              Department of Communication Engineering and Informatics, 
              The University of Electro-Communications, 1-5-1, Chofugaoka, Chofu, Tokyo, 182-8585, Japan \\
              Tel.: +81-42-443-5360 \\
              Fax: +81-42-443-5360 \\
              \email{yusaku.yamamoto@uec.ac.jp}           
           \and
           S.~Kudo \at
              Department of Communication Engineering and Informatics, 
              The University of Electro-Communications, 1-5-1, Chofugaoka, Chofu, Tokyo, 182-8585, Japan \\
               \emph{Present address: RIKEN R-CCS, 7-1-26, Minatojima-minami-machi, Chuo-ku, Kobe, Hyogo 650-0047, Japan}   
           \and
           H.~Imachi \at
              Department of Applied Mathematics and Physics, Tottori University, 4-101 Koyama-Minami, Tottori, 680-8552, Japan \\
               \emph{Present address: Preferred Networks, Inc., Otemachi Bldg., 1-6-1 Otemachi, Chiyoda-ku, Tokyo, 100-0004, Japan}   
           \and
           Y.~Miyatake \at
              Cybermedia Center, Osaka University, 1-32 Machikaneyama, Toyonaka, Osaka, 560-0043, Japan \\
           \and
           T.~Hoshi \at
              Department of Applied Mathematics and Physics, Tottori University, 4-101 Koyama-Minami, Tottori, 680-8552, Japan \\
 }

\date{Received: date / Accepted: date}

\maketitle

\begin{abstract}
In simulating physical systems, conservation of the total energy is often essential, especially when energy conversion between different forms of energy occurs frequently. Recently, a new fourth order energy-preserving integrator named MB4 was proposed based on the so-called continuous stage Runge--Kutta methods (Y.~Miyatake and J.~C.~Butcher, SIAM J.~Numer.~Anal., 54(3), 1993-2013). A salient feature of this method is that it is parallelizable, which makes its computational time for one time step comparable to that of second order methods. In this paper, we illustrate how to apply the MB4 method to a concrete ordinary differential equation using the nonlinear Schr\"{o}dinger-type equation on a two-dimensional grid as an example. This system is a prototypical model of two-dimensional disordered organic material and is difficult to solve with standard methods like the classical Runge--Kutta methods due to the nonlinearity and the $\delta$-function like potential coming from defects. Numerical tests show that the method can solve the equation stably and preserves the total energy to 16-digit accuracy throughout the simulation. It is also shown that parallelization of the method yields up to 2.8 times speedup using 3 computational nodes.
\end{abstract}

\section{Introduction}
A fundamental issue in computational physics is to develop 
parallelizable numerical algorithms for fast and reliable simulations 
and the conservation of total energy is often of vital importance 
to obtain reliable numerical results.
In particular, when the system is nonlinear and energy conversion between different types of energies, e.~g.~potential energy and kinetic energy, occurs frequently, energy conservation is critical to describe the physical process correctly. For this reason, various energy-preserving numerical schemes have been proposed both for ordinary and partial differential equations. In this paper, we focus on the numerical solution of the following system of ordinary differential equations:
\begin{eqnarray}
\frac{{\rm d}{\bm u}}{{\rm d}t} &=& {\bm f}({\bm u}), \quad {\bm u}(t_0)={\bm u}_0\in\mathbb{R}^N, \label{eq:EOM1} \\
{\bm f}({\bm u}) &=& S\nabla H({\bm u}), \label{eq:EOM2}
\end{eqnarray}
where $H({\bm u}): {\mathbb R}^N\rightarrow{\mathbb R}$ is the Hamiltonian that represents the total energy of the system and $S\in\mathbb{R}^{N\times N}$ is some constant nonsingular skew-symmetric matrix. Using \eqref{eq:EOM1} and \eqref{eq:EOM2}, it is readily verified that the total energy is preserved:
\begin{equation}
\frac{\rm d}{{\rm d}t}H({\bm u})=\nabla H\cdot\frac{d{\bm u}}{dt}=\nabla H\cdot S\nabla H({\bm u}) = 0.
\label{eq:EnergyPreservation}
\end{equation}

The simplest form of an energy-preserving method is the projection method, which projects the numerical solution computed at each time step onto a manifold with constant energy. Unfortunately, it is known that projection-based methods usually do not give satisfactory results in terms of long-range behavior~\cite{Hairer06}. A more sophisticated approach is the discrete gradient method \cite{Gonzales96,McLachlan99}, which discretizes the equation in such a way that the mechanism of energy preservation (\ref {eq:EnergyPreservation}) is maintained. The average vector field (AVF) method \cite{Quispel08}, which is of order two, belongs to this class. As an extension of the AVF method, the AVF collocation method \cite{Hairer10} has been proposed. Using the framework of the AVF collocation method, it is possible to derive higher order methods systematically. However, they are implicit schemes and the size of the system of nonlinear equations to be solved at each time step grows with the required order. This incurs large computational cost and makes it difficult to apply higher order energy-preserving schemes to large-scale problems.

Recently, Miyatake and Butcher proposed a new 4th order energy-preserving scheme that is applicable to \eqref{eq:EOM1}\eqref{eq:EOM2} based on the so-called continuous stage Runge--Kutta (CSRK) methods \cite{Miyatake16}. A distinctive feature of this method, named MB4, is that it is {\it parallelizable}. While it requires solving a system of nonlinear equations in $3N$ unknowns at each time step, the system is naturally decomposed into three independent systems of equations in $N$ unknowns each. Thus, by using three processors, one can solve the whole system in a time to solve a system in $N$ unknowns. This greatly saves the computing time, since the work required to solve a system usually grows quadratically or cubically with the number of unknowns.

The present study is inter-disciplinary between mathematics and application (simulation) research and is motivated by two aspects. The first aspect is a general need on the side of application researchers to test and evaluate numerical time integrators. So far, many numerical integrators have been proposed and they are quite different in terms of the computational cost, numerical robustness, and so on. Since such characteristics are often difficult to predict purely theoretically, a test program or a mini-application for evaluating numerical integrators is highly desired. In particular, various kinds of structure-preserving integrators such as those listed above have been developed recently. Most application researchers are not familiar with these methods and would like to compare them experimentally, so as to choose the optimal one according to their needs. The second aspect is the focus on a nonlinear time-dependent Schr\"odinger-type (NLS-type) equation. The NLS equation is one of the most famous nonlinear partial differential equations and has been investigated by many conventional methods. Thus, it is an optimal test problem to compare different numerical integrators. Moreover, the NLS-like equation has been used in quantum-mechanical wavepacket simulations and application researchers have a need for a reliable and fast numerical integrator with parallelism.


The objective of this paper is twofold. The first goal is to illustrate how to apply the MB4 method to a concrete ordinary differential equation using the nonlinear Schr\"{o}dinger-type equation as an example. Whereas the CSRK methods have several advantages over conventional Runge--Kutta methods, they require the evaluation of some definite integrals involving ${\bm f}$, which must be done either analytically or numerically. We will explain this process in detail so that an interested reader can apply the MB4 method to other ordinary differential equations as well. The second goal is to evaluate the numerical properties and parallel performance of the MB4 method. We first evaluate the accuracy and stability of the method by comparing it with the classical Runge--Kutta method (RK4), two other energy-preserving methods (AVF2 and AVF4) and two symplectic methods (GAUSS2 and GAUSS4). Besides the total energy, our nonlinear Schr\"{o}dinger-type equation has one more conserved quantity, namely, the total probability. So it is of interest to see to what extent this quantity is conserved with the MB4 method. Next, we propose an efficient implementation of the MB4 method using a sparse direct solver, parallelize it using MPI and evaluate its performance. Experimental results show that a relatively large system on a $640\times 640$ grid can be simulated in a practical time.

It is to be noted that the nonlinear Schr\"{o}dinger equation is an important differential equation with a rich mathematical structure and a wide range of applications and therefore many special-purpose numerical schemes have been designed for that equation; see \cite{Delfour81,Besse04} for examples and \cite{Sanz-Serna86} for a survey. However, our focus here is to evaluate MB4 as a general-purpose energy-preserving scheme and the equation was chosen merely as an example. Hence, we have selected only general-purpose numerical schemes that are applicable to \eqref{eq:EOM1}\eqref{eq:EOM2} for comparison.

The rest of this paper is structured as follows. In Section 2, we review the MB4 method and show how it can be applied to our nonlinear Schr\"{o}dinger-type equation. Experimental results that demonstrate the accuracy, stability and parallel performance of our implementation are given in Section 3. Finally, Section 4 concludes the paper and provides some future research directions. Throughout the paper, we use $\bar{\bm u}$, ${\bm u}^{\top}$ and ${\bm u}^{\mathsf{H}}$ to denote the complex conjugate, transpose and Hermitian conjugate of a vector ${\bm u}$, respectively. The $k$th element of ${\bm u}$ is denoted by $({\bm u})_k$ or $u_k$.

\section{Application of the MB4 method to the nonlinear Schr\"{o}dinger-type equation}
\subsection{The nonlinear Schr\"{o}dinger-type equation}
The ordinary differential equation we consider in this paper is a space-discretized version of the two-dimensional nonlinear Schr\"{o}dinger equation:
\begin{equation}
i\frac{\partial u}{\partial t}=-\nabla^2 u-\gamma|u|^2 u+vu,
\label{eq:NLS2d}
\end{equation}
where $u=u({\bm r},t)$ is the wave function, $\gamma>0$ is a scalar parameter and $v=v({\bm r})$ is an external potential. We discretize this equation on an $N_x\times N_y$ regular grid with periodic boundary conditions using the 5-point finite difference formula for the Laplacian. The resulting equation is
\begin{equation}
i\frac{{\rm d}}{{\rm d}t}{\bm u}=K{\bm u}-\gamma(\bar{\bm u}\circ{\bm u})\circ{\bm u}+V{\bm u}.
\label{eq:dNLS}
\end{equation}
Here, ${\bm u}\in\mathbb{R}^{N_x N_y}$ is a vector whose $((j-1)N_x+i)$th element is the wave function value at the $(i,j)$ grid point, $K\in\mathbb{R}^{(N_x N_y)\times(N_x N_y)}$ is an Hermitian matrix representing the kinetic energy operator defined by
\begin{equation}
K=-\left(I_{N_y}\otimes \left(\frac{1}{h_x^2}D_{N_x}\right)+\left(\frac{1}{h_y^2}D_{N_y}\right)\otimes I_{N_x}\right),
\label{eq:Kdefinition}
\end{equation}
where $h_x$ and $h_y$ are grid spacing in the $x$ and $y$ directions, respectively, $D_N\in\mathbb{R}^{N\times N}$ is a circulant matrix with $(D_N)_{i,i}=-2$ and $(D_N)_{i,i+1}=(D_N)_{i,i-1}=1$, $I_N$ is the identity matrix of order $N$, $V\in\mathbb{R}^{(N_x N_y)\times(N_x N_y)}$ is a diagonal matrix whose diagonal elements are the external potential at each grid point and $\circ$ denotes the Hadamard (componentwise) product of two vectors. Eq.~(\ref{eq:dNLS}) can be viewed as a model of two-dimensional disordered organic material, where each grid point corresponds to an organic molecule and $V$ is a potential from defects. The nonlinear term $-\gamma(\bar{\bm u}\circ{\bm u})\circ{\bm u}$ is proposed as a phenomenological one to take into account the environmental effect, since real material consists of many electrons. This term causes a lower potential value at a point where the probability density $|u|^2$ is high, and therefore has the effect of localizing the wavefunction. Eq.~(\ref{eq:dNLS}) can be expressed compactly as
\begin{equation}
i\frac{{\rm d}}{{\rm d}t}{\bm u} = \nabla_{\bar{\bm u}}H,
\label{eq:dNLS2}
\end{equation}
where
\begin{eqnarray}
H &=& U_K+U_I+U_E, \\
U_K &=& {\bm u}^{\mathsf{H}}K{\bm u}, \label{eq:U_K} \\
U_I &=& -\frac{\gamma}{2}\|{\bm n}(t)\|^2, \\
U_E &=& {\bm u}^{\mathsf{H}}V{\bm u}, \label{eq:U_E} \\
{\bm n}(t) &=& (|u_1(t)|^2,\ldots,|u_{N^2}(t)|^2)^{\top}.
\end{eqnarray}
Now, let the real and imaginary parts of ${\bm u}$ be denoted by ${\bm p}$ and ${\bm q}$, respectively. Then, (\ref{eq:dNLS2}) can be rewritten as
\begin{equation}
\frac{\rm d}{{\rm d}t}\left[\begin{array}{c}
{\bm p} \\
{\bm q}
\end{array}
\right]
=\left[\begin{array}{cc}
O & \frac{1}{2}I_{N_x N_y} \\
-\frac{1}{2}I_{N_x N_y} & O
\end{array}
\right]
\left[\begin{array}{c}
\nabla_{{\bm p}} \\
\nabla_{{\bm q}}
\end{array}
\right]H,
\end{equation}
showing that (\ref{eq:dNLS2}) is of the form \eqref{eq:EOM1}\eqref{eq:EOM2}.

In concluding this subsection, we will make a comment on the form of the kinetic energy operator $K$. While we assume that $K$ is a 5-point finite difference operator defined by (\ref{eq:Kdefinition}) in this paper, other forms of $K$ might be more appropriate in other situations. We stress that the specific form of $K$ is irrelevant in the application of the MB4 method to be described below, as long as $K$ is Hermitian. An application of this fact will be described in the Appendix.

\subsection{The MB4 method}
\label{MB4}
The MB4 method is a variant of the so-called CSRK methods that is both energy-preserving and parallelizable. In this subsection, we describe the application of this method to the system \eqref{eq:EOM1}\eqref{eq:EOM2} briefly. In the CSRK methods, the numerical solution ${\bm u}_1$ at time $t_0+h$ is computed from the solution ${\bm u}_0$ at time $t_0$ by the following formula.\begin{eqnarray}
{\bm U}_{\tau} &=& {\bm u}_0+h\int_0^1 A_{\tau,\zeta}{\bm f}({\bm U}_{\zeta})\,{\rm d}\zeta, \label{eq:MB4_1} \\
{\bm u}_1 &=& {\bm u}_0+h\int_0^1 B_{\tau}{\bm f}({\bm U}_{\tau})\,{\rm d}\tau. \label{eq:MB4_2}
\end{eqnarray}
Here, $A_{\tau,\zeta}$ is a polynomial of order $s$ in $\tau$ and of order $s-1$ in $\zeta$ satisfying $A_{0,\zeta}=0$ and $B_{\zeta}=A_{1,\zeta}$. In the MB4 method, $s=3$ and $A_{\tau,\zeta}$ is defined as
\begin{equation}
A_{\tau,\zeta}=\left[\tau\;\; \frac{\tau^2}{2}\;\; \frac{\tau^3}{3}\right]
\left[\begin{array}{ccc}
\alpha_1+4 & -6\alpha_1-6 & 6\alpha_1 \\
-6\alpha_1-6 & 36\alpha_1+12 & -36\alpha_1 \\
6\alpha_1 & -36\alpha_1 & 36\alpha_1
\end{array}
\right]
\left[\begin{array}{c}
1 \\ \zeta \\ \zeta^2
\end{array}
\right],
\label{eq:Adefinition}
\end{equation}
where $\alpha_1$ is chosen to satisfy $\alpha_1<-300\times 0.7770503941$ \cite{Miyatake16}. Hence, ${\bm U}_{\tau}$ is a vector whose components are third order polynomials in $\tau$. Now, let us express ${\bm U}_{\tau}$ using its values at four points $c_0=0, c_1, c_2, c_3$ by Lagrange interpolation:
\begin{equation}
{\bm U}_{\tau}={\bm u}_0 l_0(\tau) + \sum_{i=1}^3{\bm U}_{c_i} l_i(\tau),
\label{eq:Lagrange}
\end{equation}
where $l_i(\tau)$ is a third order polynomial satisfying $l_i(c_j)=\delta_{ij}$. Plugging (\ref {eq:Lagrange}) into (\ref{eq:MB4_1}) and noting that (\ref{eq:MB4_1}) needs to be satisfied only at $\tau=c_1, c_2, c_3$ (since ${\bm U}_{\tau}$ is a third order polynomial and (\ref{eq:MB4_1}) holds at $\tau=c_0=0$ trivially), we have the following nonlinear simultaneous equations in ${\bm U}_{c_i}$ ($i=1, 2, 3$):
\begin{equation}
{\bm U}_{c_i}={\bm u}_0+h\int_0^1 A_{c_i,\zeta}{\bm f}\left({\bm u}_0 l_0(\zeta) + \sum_{j=1}^3{\bm U}_{c_j}l_j(\zeta)\right){\rm d}\zeta, \quad i=1, 2, 3.
\label{eq:MB4_3}
\end{equation}

We solve (\ref{eq:MB4_3}) by the simplified Newton method. Let ${\bm U}=[{\bm U}_{c_1}^{\top}, {\bm U}_{c_2}^{\top}, {\bm U}_{c_3}^{\top}]^{\top}$ and define a function $\Phi: \mathbb{R}^{3N}\rightarrow\mathbb{R}^{3N}$ by
\begin{equation}
\Phi({\bm U})={\bm U}-{\bm e}_3\otimes{\bm u}_0
-h\left[\begin{array}{c}
\int_0^1 A_{c_1,\zeta}{\bm f}\left({\bm U}_{\zeta}\right){\rm d}\zeta \\
\int_0^1 A_{c_2,\zeta}{\bm f}\left({\bm U}_{\zeta}\right){\rm d}\zeta \\
\int_0^1 A_{c_3,\zeta}{\bm f}\left({\bm U}_{\zeta}\right){\rm d}\zeta
\end{array}
\right],
\label{eq:Phidefinition}
\end{equation}
where ${\bm e}_3=(1,1,1)^{\top}$. Then, (\ref{eq:MB4_3}) can be written concisely as $\Phi({\bm U})={\bf 0}$. Its Jacobian matrix is
\begin{equation}
\frac{\partial\Phi}{\partial{\bm U}}
=I_{3N}-h\left[\begin{array}{ccc}
\int_0^1 A_{c_1,\zeta}\frac{\partial{\bm f}}{\partial{\bm u}}l_1(\zeta){\rm d}\zeta  & \cdots & \int_0^1 A_{c_1,\zeta}\frac{\partial{\bm f}}{\partial{\bm u}}l_3(\zeta){\rm d}\zeta \\
\vdots & \ddots & \vdots \\
\int_0^1 A_{c_3,\zeta}\frac{\partial{\bm f}}{\partial{\bm u}}l_1(\zeta){\rm d}\zeta
& \cdots & \int_0^1 A_{c_3,\zeta}\frac{\partial{\bm f}}{\partial{\bm u}}l_3(\zeta){\rm d}\zeta\end{array}
\right].
\end{equation}
Here, $\frac{\partial{\bm f}}{\partial{\bm u}}$ in the integral is a function of ${\bm U}_{\zeta}$ and therefore of $\zeta$. In the simplified Newton method, we approximate it by its value at $\zeta=0$, namely, at ${\bm u}={\bm u}_0$, and put it outside the integral. Then, the approximated Jacobian has the following tensor product structure:
\begin{eqnarray}
\frac{\partial\Phi}{\partial{\bm U}}
&\simeq& I_{3N}-h\left[\begin{array}{ccc}
\left.\frac{\partial{\bm f}}{\partial{\bm u}}\right|_{{\bm u}={\bm u}_0}\int_0^1 A_{c_1,\zeta}l_1(\zeta){\rm d}\zeta  & \cdots & \left.\frac{\partial{\bm f}}{\partial{\bm u}}\right|_{{\bm u}={\bm u}_0}\int_0^1 A_{c_1,\zeta}l_3(\zeta){\rm d}\zeta \\
\vdots & \ddots & \vdots \\
\left.\frac{\partial{\bm f}}{\partial{\bm u}}\right|_{{\bm u}={\bm u}_0}\int_0^1 A_{c_3,\zeta}l_1(\zeta){\rm d}\zeta
& \cdots & \left.\frac{\partial{\bm f}}{\partial{\bm u}}\right|_{{\bm u}={\bm u}_0}\int_0^1 A_{c_3,\zeta}l_3(\zeta){\rm d}\zeta
\end{array}
\right] \nonumber \\
&=& I_{3N}-hE\otimes J,
\label{eq:Jacobian2}
\end{eqnarray}
where $E$ is a $3\times 3$ matrix defined as
\begin{equation}
E_{ij}=\int_0^1 A_{c_i,\zeta}l_j(\zeta){\rm d}\zeta
\label{eq:defE}
\end{equation}
and
\begin{equation}
J=\left.\frac{\partial{\bm f}}{\partial{\bm u}}\right|_{{\bm u}={\bm u}_0}.
\end{equation}
Using the approximated Jacobian of (\ref{eq:Jacobian2}), the iteration of the simplified Newton method can be written as
\begin{equation}
\frac{\partial\Phi}{\partial{\bm U}}\,{\bm r}^l=-\Phi({\bm U}^l), \quad {\bm U}^{l+1}={\bm U}^l+{\bm r}^l, \quad l=0, 1, 2,\ldots.
\label{eq:Newton1}
\end{equation}

In the MB4 method, the matrix $E$ is diagonalizable and has real eigenvalues that do not depend on the choice of $c_1, c_2, c_3$ (\cite[Theorem 5.1]{Miyatake16}). Thus, using the eigendecomposition $E=T\Lambda T^{-1}$, where $T\in\mathbb{R}^{3\times 3}$ is nonsingular and $\Lambda\in\mathbb{R}^{3\times 3}$ is diagonal, the first equation of (\ref{eq:Newton1}) can be written as
\begin{equation}
\left\{I_{3N}-h(T\Lambda T^{-1})\otimes J\right\}{\bm r}^l=-\Phi({\bm U}^l),
\label{eq:Newton2}
\end{equation}
or
\begin{equation}
(I_{3N}-h\Lambda\otimes J)\left\{(T^{-1}\otimes I_N){\bm r}^l\right\}=-(T^{-1}\otimes I_N)\Phi({\bm U}^l).
\label{eq:Newton3}
\end{equation}
Since the coefficient matrix $I_{3N}-h\Lambda\otimes J$ of the linear system (\ref {eq:Newton3}) has a $3\times 3$ block diagonal structure, the system can be solved in the following three steps.
\begin{eqnarray}
&& \bar{\Phi}^l=(T^{-1}\otimes I_N)\Phi({\bm U}^l), \label{eq:MB4_4} \\
&& \bar{\bm r}_i^l=-(I_N-h\lambda_i J)^{-1}\bar{\Phi}_i^l, \quad i=1, 2, 3, \label{eq:MB4_5} \\
&& {\bm r}^l=(T\otimes I_N)\bar{\bm r}^l.
\label{eq:MB4_6}
\end{eqnarray}
Here, $\bar{\Phi}_i^l$ and $\bar{\bm r}_i^l$ are $N$-dimensional vectors consisting of the $((i-1)N+1)$th through the $(iN)$th elements of $\bar{\Phi}^l$ and $\bar{\bm r}^l$, respectively. The solution of Eq.~(\ref{eq:MB4_5}), which accounts for most of the computational work, amounts to solving three independent linear systems of order $N$ each and therefore can be performed by three processors in parallel.

It is to be noted that the simplified Newton method can fail to converge if $h$ is too large. This can be detected from the fact that $\|{\bm r}^l\|$ does not become sufficiently small even after a pre-specified number of iterations. If this occurs, it is advisable to decrease $h$ and redo the simplified Newton iteration.

Once ${\bm U}=[{\bm U}_{c_1}^{\top}, {\bm U}_{c_2}^{\top}, {\bm U}_{c_3}^{\top}]^{\top}$ has been obtained, we can recover ${\bm U}_{\tau}$ by inserting ${\bm U}_{c_1}, {\bm U}_{c_2}, {\bm U}_{c_3}$ into (\ref{eq:Lagrange}) and then use (\ref{eq:MB4_2}) to get ${\bm u}_1$. However, when $c_3=1$, we can omit these steps since we have
\begin{equation}
{\bm u}_1 = {\bm u}_0+h\int_0^1 B_{\tau}{\bm f}({\bm U}_{\tau})\,{\rm d}\tau = {\bm u}_0+h\int_0^1 A_{1,\tau}{\bm f}({\bm U}_{\tau})\,{\rm d}\tau
= {\bm U}_{c_3}
\end{equation}
from $B_{\tau}=A_{1,\tau}$ and the definition of ${\bm U}_{c_3}$. Hence, ${\bm U}_{c_3}$ can be used directly as the approximate solution at the next time step.

\subsection{The MB4 method applied to the nonlinear Schr\"{o}dinger-type equation}
Now we apply the MB4 method to the nonlinear Schr\"{o}dinger-type equation (\ref {eq:dNLS}). We first rewrite (\ref {eq:dNLS}) as
\begin{eqnarray}
\frac{{\rm d}}{{\rm d}t}{\bm u} &=& {\bm f}_K({\bm u})+{\bm f}_I({\bm u})+{\bm f}_E({\bm u}), \label{eq:fdivide} \\
{\bm f}_K({\bm u}) &=& -iK{\bm u}, \label{eq:fKu} \\
{\bm f}_I({\bm u}) &=& i\gamma(\bar{\bm u}\circ{\bm u})\circ{\bm u} \label{eq:fIu} \\
{\bm f}_E({\bm u}) &=& -iV{\bm u} \label{eq:fEu}.
\end{eqnarray}
Now we define ${\bm U}_{\tau}, {\bm U}_{c_1}, {\bm U}_{c_2}, {\bm U}_{c_3}, {\bm u}_0, {\bm u}_1$ and ${\bm U}$ as we did in the previous subsection and solve the equation $\Phi({\bm U})={\bm 0}$ corresponding to (\ref {eq:MB4_3}). To this end, we need to derive formulas to compute $\Phi({\bm U})$ and $J=\left.\partial{\bm f}/\partial{\bm u}\right|_{{\bm u}={\bm u}_0}$ as functions of ${\bm U}$. Their derivation will be explained in the following.

\paragraph{Computation of $\Phi({\bm U})$}
To derive the concrete form of $\Phi({\bm U})$, we need to compute $\int_0^1 A_{c_i,\zeta}{\bm f} ({\bm U}_\zeta)\,{\rm d}\zeta$ for $i=1, 2, 3$ and ${\bm f}={\bm f}_K, {\bm f}_I$ and ${\bm f}_E$. Note that ${\bm f}_K({\bm u})$ is a term that is linear in ${\bm u}$, while ${\bm f}_I({\bm u})$ is a term that is nonlinear in ${\bm u}$. We first consider the linear terms ${\bm f}_K({\bm u})$ and ${\bm f}_E({\bm u})$ together:
\begin{eqnarray}
&& \int_0^1 A_{c_i,\zeta}\left\{{\bm f}_K({\bm U}_\zeta)+{\bm f}_E({\bm U}_\zeta)\right\}\,{\rm d}\zeta \nonumber \\
&&= -i\int_0^1 A_{c_i,\zeta}(K+V){\bm U}_\zeta\,{\rm d}\zeta \nonumber \\
&&= -i(K+V)\int_0^1 A_{c_i,\zeta}{\bm U}_\zeta\,{\rm d}\zeta \nonumber \\
&&= -i(K+V)\int_0^1 A_{c_i,\zeta}\left({\bm u}_0 l_0(\zeta) + \sum_{j=1}^3{\bm U}_{c_j} l_j(\zeta)\right)\,{\rm d}\zeta \nonumber \\
&&= -i(K+V)
\left\{
{\bm u}_0\int_0^1 A_{c_i,\zeta}l_0(\zeta)\,{\rm d}\zeta + \sum_{j=1}^3{\bm U}_{c_j}\int_0^1 A_{c_i,\zeta}l_j(\zeta)\,{\rm d}\zeta
\right\} \nonumber \\
&&= -i(K+V)\left(E_{i0}{\bm u}_0+\sum_{j=1}^3 E_{ij}{\bm U}_{c_j}\right) \q (i=0, 1, 2, 3),
\label{eq:fKu2}
\end{eqnarray}
where $E_{ij}$ ($1\le i,j\le 3$) is given by (\ref{eq:defE}) and
\begin{equation}
E_{i0} = \int_0^1 A_{c_i,\zeta}l_0(\zeta)\,{\rm d}\zeta \q (i=0, 1, 2, 3).
\label{eq:Ei0}
\end{equation}
Since $E_{ij}$'s of Eqs.~(\ref{eq:defE}) and (\ref{eq:Ei0}) are constants that depend only on $\alpha_1$ and $c_i$'s (which are constants throughout the whole computation), they can be pre-computed. Thus, all we need is to compute the linear combination $E_{i0}{\bm u}_0+\sum_{j=1}^3 E_{ij}{\bm U}_{c_j}$ and multiply it by the matrix $-i(K+V)$.

On the other hand, ${\bm f}_I({\bm u})$ can be computed independently for each grid point and its value at the $k$th grid point ($1\le k\le N_x N_y$) is $({\bm f}_I({\bm u}))_k=i\gamma|u_k|^2 u_k$. Hence, the $k$th component of $\int_0^1 A_{c_i,\zeta}{\bm f}_I ({\bm U}_\zeta)\,{\rm d}\zeta$ is computed as
\begin{eqnarray}
&& \int_0^1 A_{c_i,\zeta}({\bm f}_I({\bm U}_\zeta))_k\,{\rm d}\zeta \nonumber \\
&&= i\gamma\int_0^1 A_{c_i,\zeta}|U_{\zeta,k}|^2 U_{\zeta,k}\,{\rm d}\zeta \nonumber \\
&&= i\gamma\int_0^1 A_{c_i,\zeta}\left|u_{0,k} l_0(\zeta) + \sum_{i=1}^3U_{c_i,k} l_i(\zeta)\right|^2 \left(u_{0,k} l_0(\zeta) + \sum_{i=1}^3 U_{c_i,k} l_i(\zeta)\right)\,{\rm d}\zeta. \nonumber \\
\label{eq:fIu2}
\end{eqnarray}
This is a third order homogeneous polynomial in $u_{0,k}, U_{c_1,k}, U_{c_2,k}, U_{c_3,k}$ and their complex conjugate. If we regard the real and imaginary parts of these variables as separate variables, we obtain a third order homogeneous polynomial in eight variables. Its coefficients are given as definite integrals of 11th order polynomials in $\zeta$ (since $A_{c_i,\zeta}$ and $l_i(\zeta)$ are second and third order polynomials in $\zeta$, respectively). Because these coefficients are constants that do not depend on $k$, they can be pre-computed.

\paragraph{Computation of $J$}
We write the Jacobian $J$ as
\begin{equation}
J=\left.\frac{\partial{\bm f}}{\partial{\bm u}}\right|_{{\bm u}={\bm u}_0}=J_K+J_I+J_E,
\label{eq:Jacobianall}
\end{equation}
where
\begin{equation}
J_K=\left.\frac{\partial{\bm f}_K}{\partial{\bm u}}\right|_{{\bm u}={\bm u}_0}, \q
J_I=\left.\frac{\partial{\bm f}_I}{\partial{\bm u}}\right|_{{\bm u}={\bm u}_0}, \q
J_E=\left.\frac{\partial{\bm f}_E}{\partial{\bm u}}\right|_{{\bm u}={\bm u}_0}.
\label{eq:Jacobian} 
\end{equation}
$J_K$ and $J_E$ can be computed from (\ref{eq:fKu}) and (\ref{eq:fEu}) as
\begin{equation}
J_K+J_E=-i(K+V).
\end{equation}
If we deal with the real and imaginary parts as separate variables and write ${\bm u}={\bm p}+i{\bm q}$, ${\bm f}_K={\bm g}_K+i{\bm h}_K$ and ${\bm f}_E={\bm g}_E+i{\bm h}_E$, we have
\begin{equation}
J_K+J_E = \left[\begin{array}{cc}
\frac{\partial{\bm g}_K}{\partial{\bm p}}+\frac{\partial{\bm g}_E}{\partial{\bm p}} &
\frac{\partial{\bm g}_K}{\partial{\bm q}}+\frac{\partial{\bm g}_E}{\partial{\bm q}} \\
\frac{\partial{\bm h}_K}{\partial{\bm p}}+\frac{\partial{\bm h}_E}{\partial{\bm p}}
&
\frac{\partial{\bm h}_K}{\partial{\bm q}}+\frac{\partial{\bm h}_E}{\partial{\bm q}}
\end{array}
\right]_{{\bm p}={\bm p}_0, {\bm q}={\bm q}_0}
=\left[\begin{array}{cc}
O & A \\
-A & O
\end{array}
\right],
\label{eq:JL}
\end{equation}
where $A=K+V$. Next, we compute $J_I$. By writing ${\bm f}_I={\bm g}_I+i{\bm h}_I$, we obtain
\begin{eqnarray}
({\bm g}_I)_k &=& -\gamma(p_k^2 q_k+q_k^3), \nonumber \\
({\bm h}_I)_k &=& \gamma(p_k^3 + p_k q_k^2).
\end{eqnarray}
Hence, $J_I$ can be computed as
\begin{equation}
J_I = \left[\begin{array}{cc}
\frac{\partial{\bm g}_I}{\partial{\bm p}} & \frac{\partial{\bm g}_I}{\partial{\bm q}} \\
\frac{\partial{\bm h}_I}{\partial{\bm p}} & \frac{\partial{\bm h}_I}{\partial{\bm q}}
\end{array}
\right]_{{\bm p}={\bm p}_0, {\bm q}={\bm q}_0}
= \gamma\left[\begin{array}{cc}
{\rm diag}(-2p_k q_k) & {\rm diag}(-p_k^2-3q_k^2)\\
{\rm diag}(3p_k^2+q_k^2) & {\rm diag}(2p_k q_k)
\end{array}
\right]_{{\bm p}={\bm p}_0, {\bm q}={\bm q}_0}
\label{eq:JR}
\end{equation}
Here, ${\rm diag}(2p_k q_k) \equiv {\rm diag}(2p_1 q_1, \ldots, 2p_{N_x N_y}q_{N_x N_y})$ and so on.

Finally, $J$ is obtained by inserting (\ref{eq:JL}) and (\ref{eq:JR}) into (\ref{eq:Jacobianall}).

\vspace{3mm}

The algorithm of the MB4 method applied to the nonlinear Schr\"{o}dinger-type equation is summarized as Algorithm \ref{MB4_algorithm}. Here, it is assumed that $c_3=1$.

\begin{algorithm}[H]
\caption{MB4 method for the nonlinear Shr\"{o}dinger-type equation}\label{MB4_algorithm}
\begin{algorithmic}[1]
\State Compute $E_{ij}$ ($0\le i\le 3$, $1\le j\le 3$) by (\ref{eq:defE}) and (\ref{eq:Ei0}).
\State Compute the eigendecomposition of $E=(E_{ij})_{i,j=1}^3$: $E=T\Lambda T^{-1}$.
\State Compute the coefficients of the third order homogeneous polynomial (\ref{eq:fIu2}).
\State Set the initial wavefunction vector ${\bm u}(0)$.
\State $t:=0$
\For{$step=1, \ldots, n_{step}$} \Comment{Time evolution loop}
   \State ${\bm u}_0={\bm u}(t)$
   \State Compute the Jacobian $J$ by (\ref{eq:Jacobianall}), (\ref{eq:JL}) and (\ref{eq:JR}).
   \State ${\bm U}^0:=[{\bm u}_0^{\top}, {\bm u}_0^{\top}, {\bm u}_0^{\top}]^{\top}$. \Comment{Initial value for the simplified Newton iteration}
   \State $r=1.0$ \Comment{Measure for convergence of the simplified Newton method}
   \State $l=0$
   \While{$r>\epsilon$} \Comment{Simplified Newton loop}
      \State Compute $\Phi({\bm U}^l)$ by (\ref{eq:Phidefinition}), (\ref{eq:fKu2}) and (\ref{eq:fIu2}).
      \State $\bar{\Phi}^l:=(T^{-1}\otimes I_N)\Phi({\bm U}^l)$
      \For{$i=1, 3$} \Comment{Parallelizable loop}
         \State Solve $(I_N-h\lambda_i J)\bar{\bm r}_i^l=-\bar{\Phi}_i^l$.
      \EndFor
      \State ${\bm r}^l:=(T\otimes I_N)\bar{\bm r}^l$
      \State ${\bm U}^{l+1}:={\bm U}^l+{\bm r}^l$
      \State $r:=\|{\bm r}^l\|_{\infty}$
      \State $l:=l+1$
   \EndWhile
   \State ${\bm u}(t+h):={\bm U}_{c_3}$
   \State $t:=t+h$
\EndFor
\end{algorithmic}
\end{algorithm}

\section{Numerical experiments and software distribution}
In this section, we first evaluate the accuracy and stability of the MB4 method applied to the two-dimensional nonlinear Schr\"{o}dinger-type equation (\ref{eq:dNLS}) by comparing it with other numerical methods. Then we evaluate the parallel performance of the MB4 method on larger scale problems. Throughout this section, we use the following computational conditions unless otherwise specified.
\begin{itemize}
\item Computational domain: $(x,y)\in[0,2\pi]^2$ with periodic boundary conditions.
\item Time domain: $t\in[0, 1]$ with time step $h=0.01$.
\item Initial condition: $u(x,y)=1+2\cos x+2\cos y$.
\item Strength of the nonlinear term: $\gamma=0.1$ (see Eq.~(\ref{eq:fIu}))
\item $\alpha_1=-300\cdot\frac{19}{8}$ (see Eq.~(\ref{eq:Adefinition})) and $c_1=\frac{1}{3}, c_2=\frac{2}{3}, c_3=1$ (see subsection \ref{MB4}).
\item Convergence criterion for the simplified Newton method: $\epsilon=10^{-13}$.
\end{itemize}
The coefficients of the third order homogeneous polynomial (\ref{eq:fIu2}) were computed analytically using a formula processing system Maxima.

\subsection{Comparison with other methods \label{SEC-ACCU-STA}}
We implemented three classes of numerical methods in MATLAB using double precision arithmetic. The first class consists of energy-preserving integrators and includes the MB4, AVF2 \cite{Quispel08} and AVF4 \cite{Hairer10} methods. The second class consists of symplectic integrators and includes the GAUSS2 \cite{Hairer06} and GAUSS4 \cite{Hairer06} methods. The last class consists of only one integrator, the classical Runge--Kutta (RK4) method. Among them, the MB4, AVF4, GAUSS4 and RK4 methods are of order 4, while AVF2 and GAUSS2 are of order 2. In general, symplectic Runge--Kutta integrators preserve second order conserved quantities exactly. Thus, in our case, GAUSS2 and GAUSS4 preserve the total probability $\|{\bm u}(t)\|_2^2$, but not the total energy. On the other hand, MB4, AVF2 and AVF4 preserve the total energy exactly, but not the total probability.

\begin{figure}[h]
\centerline{\includegraphics[height=2.2in]{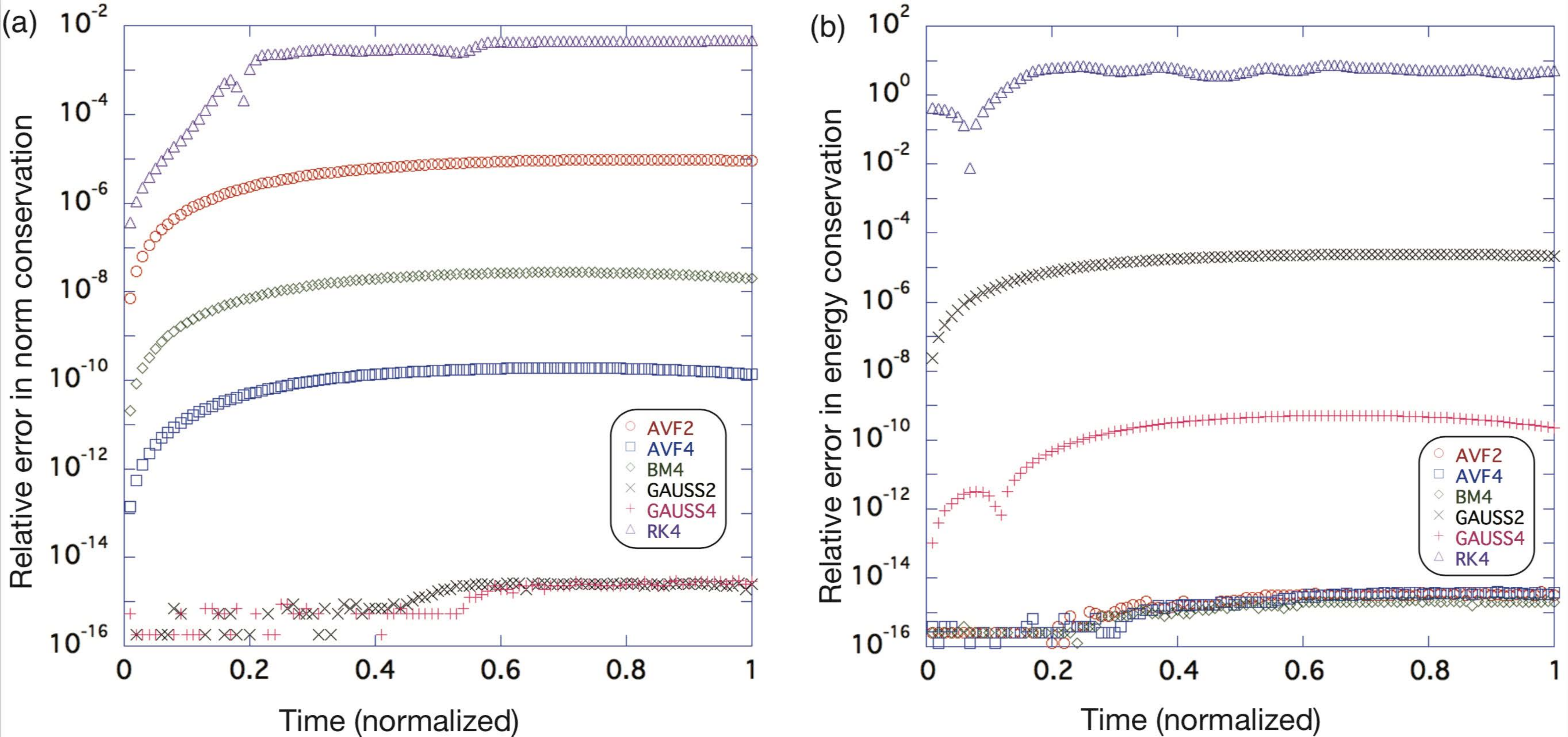}}
\caption{Deviation of the total probability (left) and the total energy (right) from the initial values.}
\label{fig:fig1}
\end{figure}

The deviation of the total probability and the total energy from their initial values is shown in Fig.~\ref{fig:fig1}(a) and \ref{fig:fig1}(b), respectively, for each numerical method. Here, the grid size is $70 \times 70$ and the potential $V(x,y)$ is set to zero. As predicted by the theory, GAUSS2 and GAUSS4 preserve the total probability and MB4, AVF2, and AVF4 preserve the total energy both to $10^{-14}$ accuracy. It is also noteworthy that MB4 and AVF4 preserves the total probability and GAUSS4 preserves the total energy both to at least $10^{-8}$ accuracy. Thus we can conclude that these fourth order structure-preserving methods can solve the nonlinear Schr\"{o}dinger-type equation both stably and accurately. In contrast, RK4 fails to preserve both total probability and total energy; their deviation from the initial values is several orders of magnitude larger than that of the fourth order structure-preserving methods.

All of the above methods except for RK4 are implicit methods and require the solution of linear simultaneous equations with a coefficient matrix of the form $I-hE\otimes J$ (see Eq.~(\ref{eq:Jacobian2})) if the simplified Newton method is used. When the number of variables is $N$, the matrix size is $N\times N$ for AVF2 and GAUSS2, $2N\times 2N$ for AVF4 and GAUSS4 and $3N\times 3N$ for MB4. However, as we have shown in subsection \ref{MB4}, the $3N\times 3N$ coefficient matrix of MB4 can be decomposed into a direct sum of three $N\times N$ matrices. Hence, theoretically, the solution of the linear simultaneous equation requires only three times the work of that for AVF2 and GAUSS2. Note that the coefficient matrices of AVF4 and GAUSS4 do not admit such decomposition, since they have complex eigenvalues.  Accordingly, these methods require $2^2 \sim 2^3$ times larger work than the AVF2 and GAUSS2 methods, depending on the complexity of linear equation solution.

\begin{figure}[h]
\centerline{\includegraphics[height=2.5in]{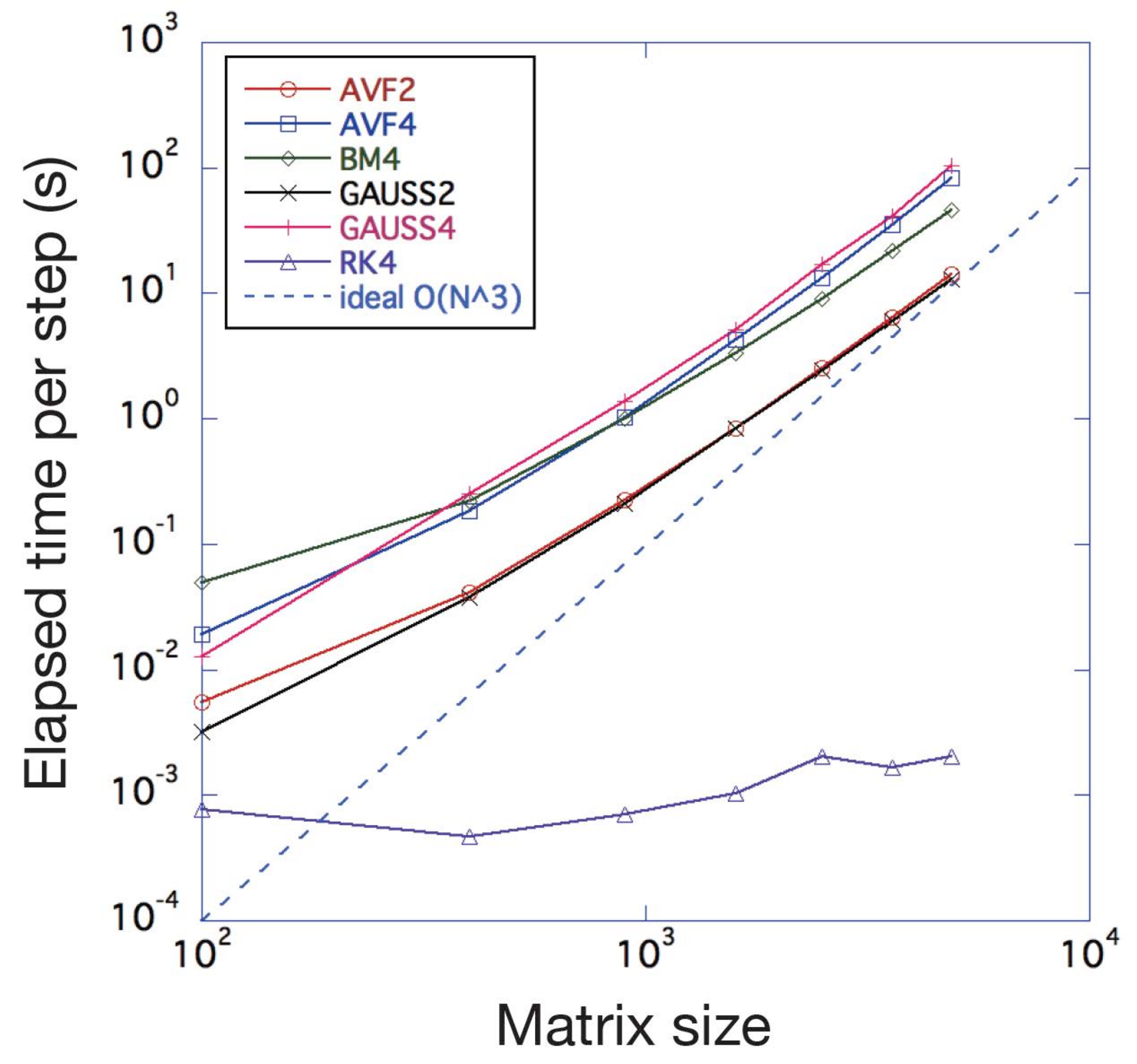}}
\caption{Computation time per step of each method (MATLAB).}
\label{fig:fig2}
\end{figure}

In Fig.~\ref{fig:fig2}, we show the computation time per step of each method as a function of the number of grid points. While the fourth order methods require more time than the second order methods, MB4 is the fastest among the former due to the reason described above. Table \ref{tab:tab1} lists the computation time per step for $70\times 70$ grid. The computation time relative to those of the second order methods is about 3 for MB4 and around 8 for AVF4 and GAUSS4. This is also consistent with the estimation given above.

\begin{table}
  \centering
  \caption{Computation time per step for $70\times 70$ grid (MATLAB).}
  \begin{tabular}{c||c|c|c|c|c|c}
    \hline
    Method & GAUSS2 & AVF2 & GAUSS4 & AVF4 & MB4 & RK4 \\ \hline\hline
    Order & 2 & 2 & 8 & 8 & 8 & 8 \\
    Time (sec) & 13.1 & 14.2 & 103.7 & 82.4 & 46.2 & 0.002 \\
    Relative time & 1 & 1 & 8 & 8 & 3 & $-$ \\ \hline
  \end{tabular}
  \label{tab:tab1}
\end{table}

\subsection{Parallel performance \label{SEC-PARALLEL}}
To evaluate parallel performance, we re-implemented the MB4 method in FORTRAN and parallelized it using MPI. In our implementation, only the {\it for} loop of lines 15--17 of Algorithm \ref{MB4_algorithm} (solution of the linear simultaneous equations) is parallelized with three MPI nodes and other parts are executed redundantly on all of these nodes to minimize inter-node data transfer. Since the solution of the linear simultaneous equations is the most computationally intensive part, even this simple parallelization scheme is expected to deliver nearly linear speedup. As a linear equation solver to run on each node, we used the PARDISO sparse direct solver included in Intel Math Kernel Library, which is multithreaded within each node. For the numerical experiments in this subsection, we used three nodes of the Reedbush-U supercomputer at the Information Technology Center of The University of Tokyo. Each node of Reedbush-U consists of two Intel Xeon E5-2695v4 (Broadwell-EP) processors, each of which has 18 cores and runs at 2.1GHz. Hence, the linear equation solution part was parallelized with $(18\times 2)$ (threads/node) $\times 3$ (nodes) $=108$ threads in total. The grid size was varied from $40\times 40$ to $640\times 640$.

\begin{figure}[h]
\centerline{\includegraphics[height=2.0in]{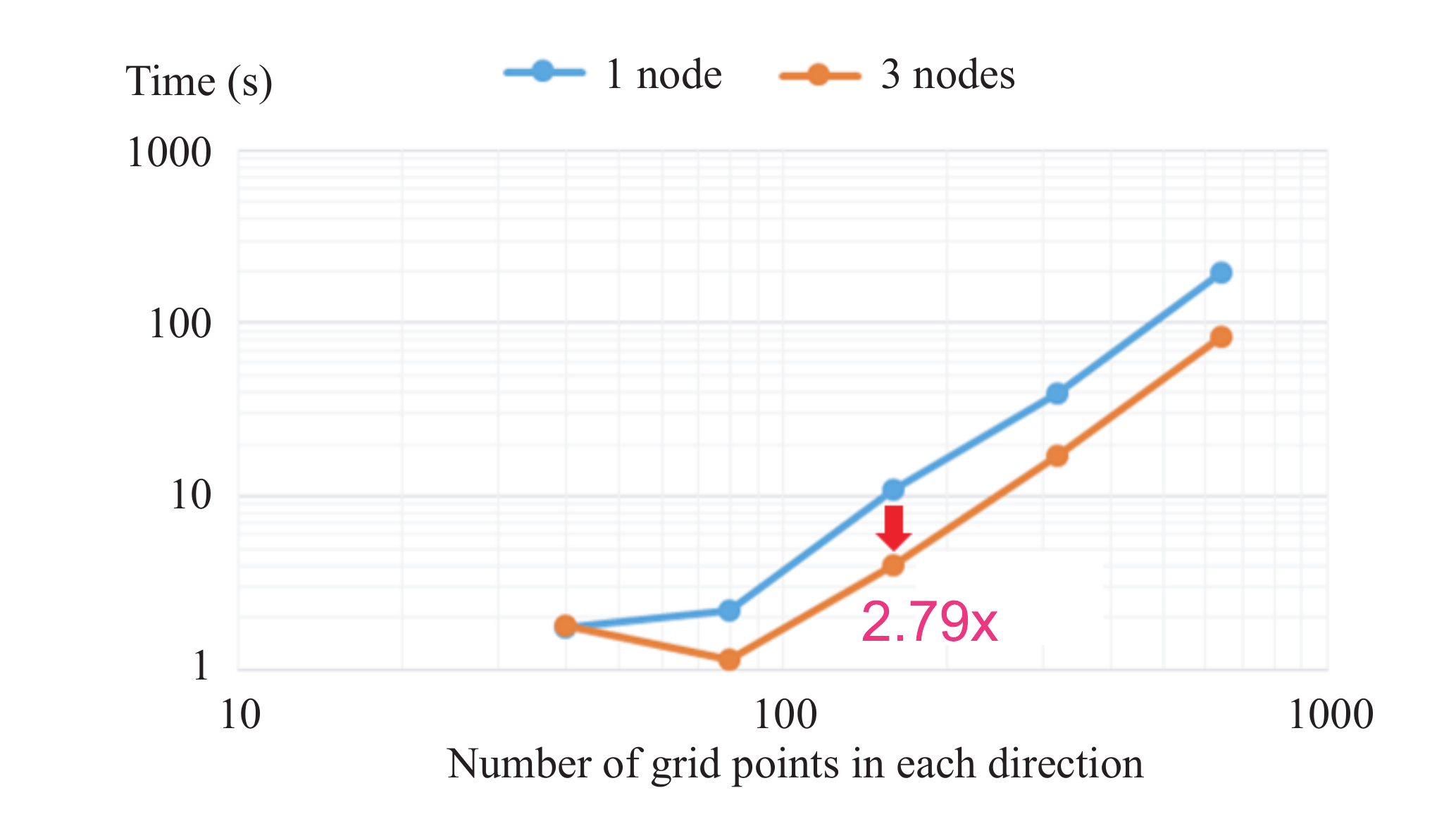}}
\caption{Computation time per step of the parallelized MB4 method (Reedbush-U).}
\label{fig:fig3}
\end{figure}

\begin{figure}[h]
\centerline{\includegraphics[height=2.0in]{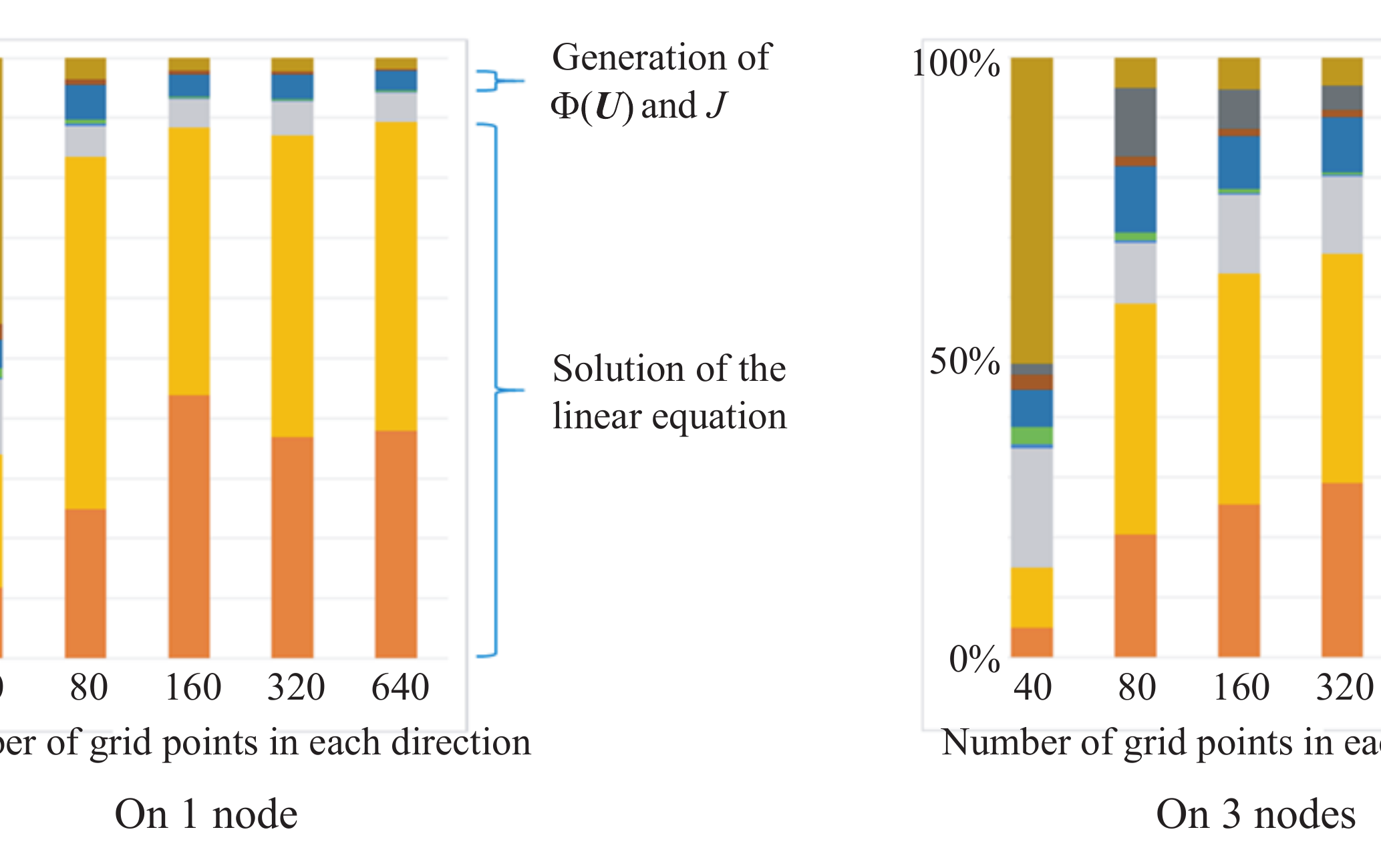}}
\caption{Breakdown of the computation time (Reedbush-U).}
\label{fig:fig4}
\end{figure}

The computation time per step is plotted in Fig.~\ref{fig:fig3} as a function of the number of grid points both for the parallel (using 3 nodes) and sequential execution. The parallel speedup in the best case is 2.79. This suggests that with the help of parallel processing, the MB4 method can provide fourth order accuracy, as well as exact energy preservation, with roughly the same computation time as that of the second order methods. The breakdown of the execution time is shown in Fig.~\ref{fig:fig4}. As expected,  nearly 90\% of the computation time is spent for the linear equation solution when the grid size is larger than $80\times 80$.

\begin{figure}[h]
\centerline{\includegraphics[height=1.0in]{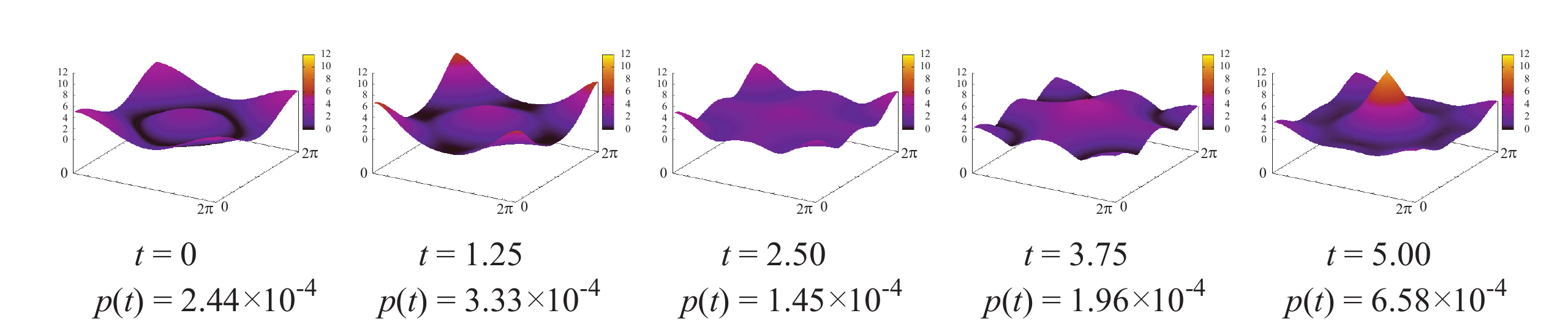}}
\caption{Plot of $|{\bm u}({\bm r},t)|^2$ for the case of $\delta$-function like potential.}
\label{fig:fig5}
\end{figure}

Finally, we made an experiment with $\delta$-function like potential, that is, $V=V_0\; (<0)$ at the origin and $V=0$ at other grid points. Physically, this can be regarded as a model of a defect or an impurity atom at the origin. In the experiment, we set the parameters to $\gamma=0.05$ and $V_0=-50$ and computed the solution up to $t=20$ with step size $h=0.01$ using a $100\times 100$ grid. The probability distribution function (PDF) $|u({\bm r},t)|^2$ at $t=1.25j$ ($0\le j\le 4$) is shown in Fig.~\ref{fig:fig5}. Due to the negative potential, the particle tends to concentrate around the origin and the PDF has a sharp peak there. Although this type of nearly discontinuous solution is difficult to deal with by numerical methods, the MB4 method can reproduce it without difficulty. The participation ratio $p(t)$, defined as the summation of $|u({\bm r},t)|^4$ over all grid points, is also shown for each $t$. The ratio is equal to $1/(N_x N_y)=10^{-4}$ when the probability density $|u({\bm r},t)|^2$ is uniform and increases as it is localized.

The temporal variations of the total energy and each component of the energy ($U_K$, $U_I$ and $U_E$, see Eqs.~(\ref{eq:U_K}) through (\ref{eq:U_E})) are shown in Fig.~6(a) and (b), respectively. Even though each component of the energy oscillates sharply due to energy conversion, the total energy is preserved to $10^{-16}$ accuracy throughout the simulation.

\begin{figure}[h]
  \begin{center}
    \begin{tabular}{c}
      \begin{minipage}{0.5\hsize}
        \begin{center}
          \includegraphics[clip, width=6cm]{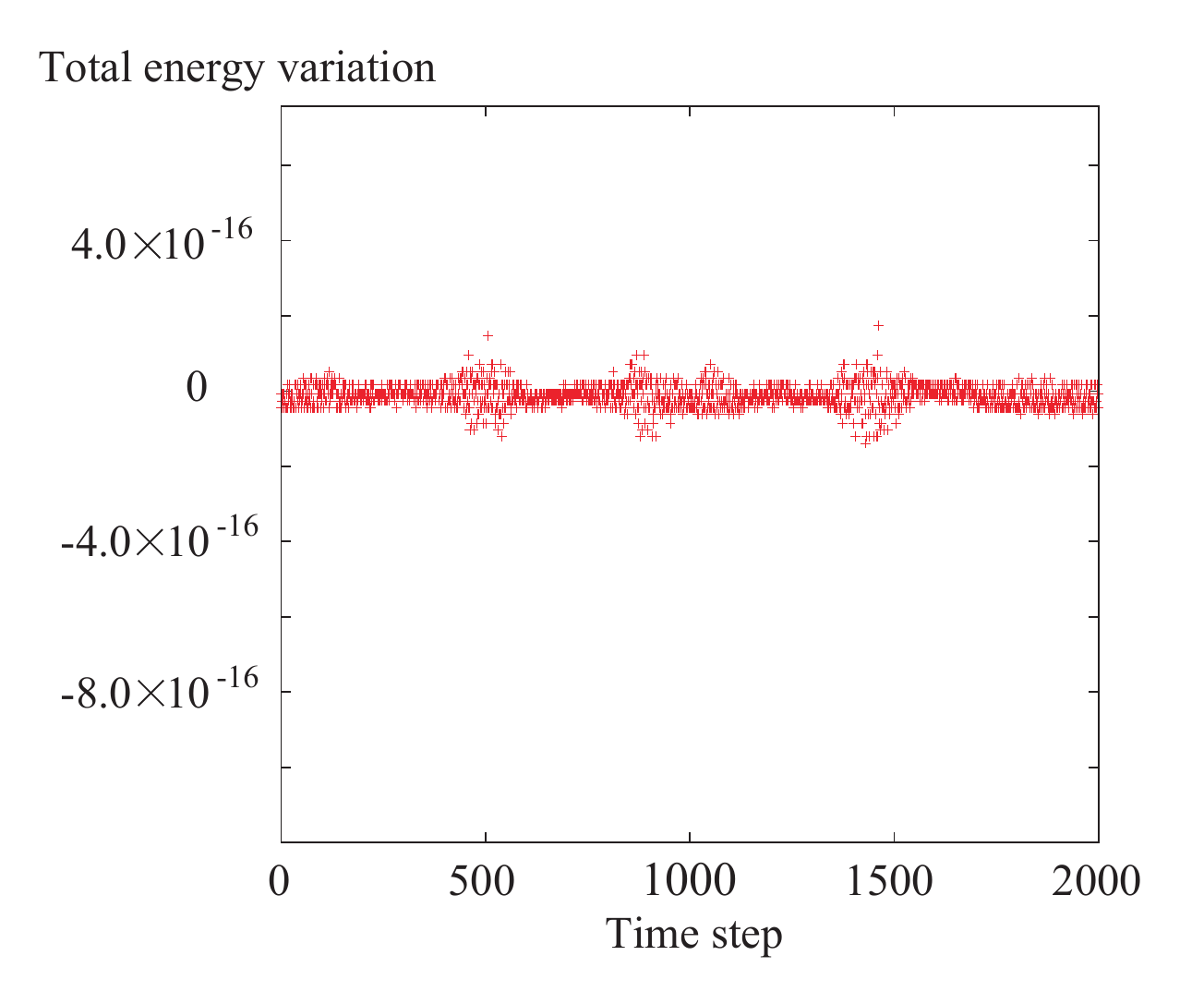}
          \hspace{1.6cm} (a) Total energy
        \end{center}
      \end{minipage}
      \begin{minipage}{0.5\hsize}
        \begin{center}
          \includegraphics[clip, width=6cm]{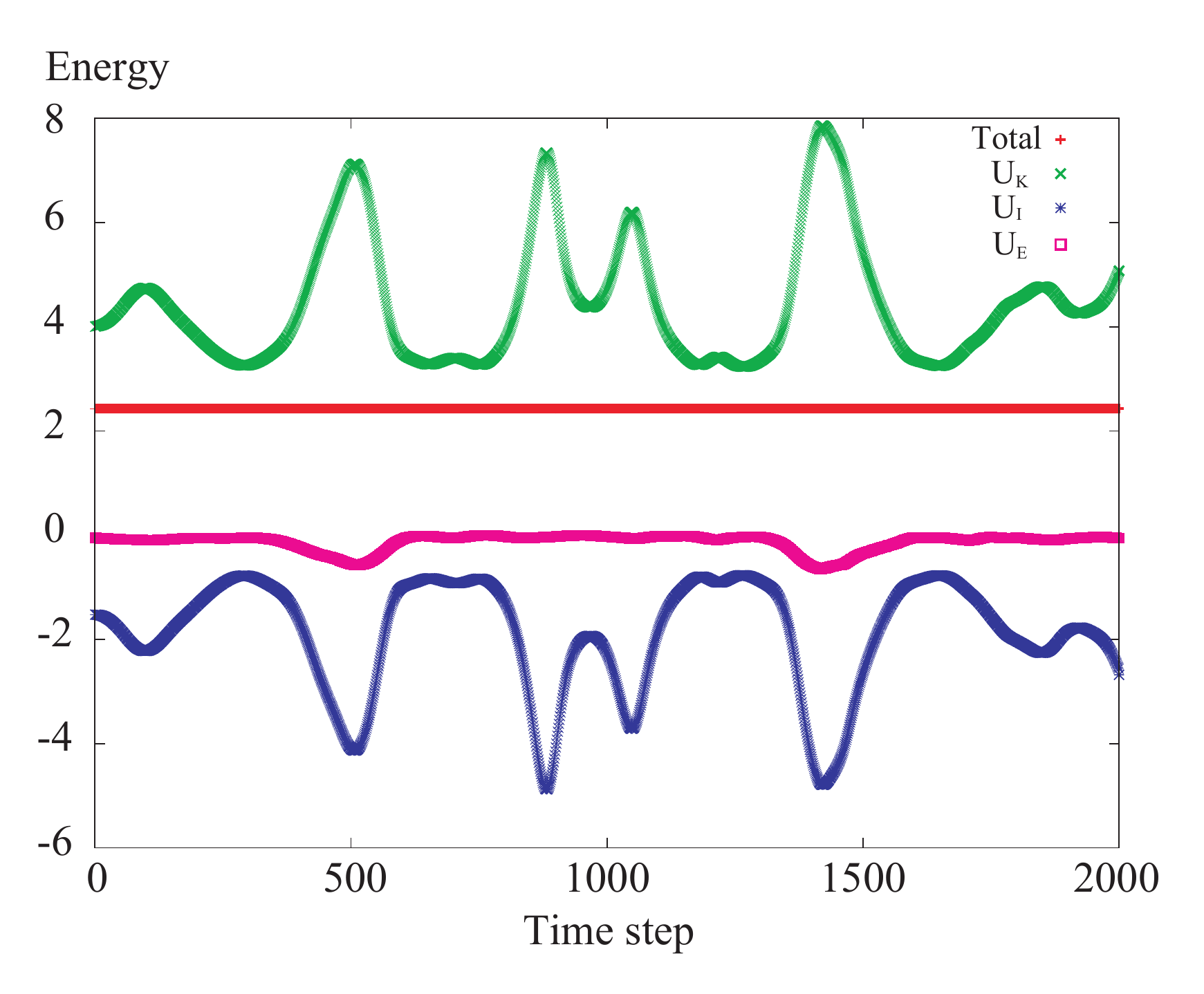}
          \hspace{1.6cm} (b) $U_K$, $U_I$ and $U_E$
        \end{center}
      \end{minipage}
    \end{tabular}
    \caption{Time evolution of the energy for the case of $\delta$-function like potential.}
  \end{center}
  \label{fig:fig6}
\end{figure}

The MB4-based parallel NLS equation solver introduced in this subsection is named {\tt nls2d\_bm4sp.f90}. The program has several parameters as listed in Table \ref{tab:parameters}. All of these parameters are defined in the main routine. By rewriting the program slightly, it is also possible to use an external potential of arbitrary functional form or a kinetic energy operator that is different from the standard one defined by Eq.~(\ref{eq:Kdefinition}).
\begin{table}[h]
  \centering
    \caption{Parameters in the program {\tt nls2d\_bm4sp.f90}.}
    \begin{tabular}{c||c|c}
      \hline
      parameter & meaning & equivalent parameter \\
      name & & in the paper \\ \hline\hline
      lx & domain size in the $x$ direction & (fixed to $2\pi$) \\
      ly & domain size in the $y$ direction & (fixed to $2\pi$) \\
      nx & grid size in the $x$ direction & $N_x$ \\
      ny & grid size in the $y$ direction & $N_y$ \\
      t\_end & simulation end time & (fixed to 1 or 20) \\
      dt & time step & $h$ \\
      eps & strength of the nonlinear term & $-\gamma$ \\
      V0 & external potential at the origin & $V_0$ \\\hline
  \end{tabular}
  \label{tab:parameters}
\end{table}

\subsection{Software distribution}

The codes used in this section
are distributed as open-source stand-alone softwares 
under the MIT license \cite{github}.
The MATLAB code used in Sec.~\ref{SEC-ACCU-STA}
is available as \lq mb4-nls2d-matlab' and
the fortran code used in Sec.~\ref{SEC-PARALLEL}
is available as \lq mb4-nls2d'. 
The codes help researchers to understand the algorithm and 
to develop their own code for different problems.

\section{Conclusion}
In this paper, we applied the MB4 method, a recently proposed 4th order numerical integrator that is both energy-preserving and parallelizable, to a nonlinear Schr\"{o}dinger-type equation on a 2-dimensional regular grid. Numerical results show that the method can solve the equation stably and accurately even in the presence of $\delta$-function like potential. Furthermore, it has been shown that an MPI version of the solver can achieve nearly linear parallel speedup and provides 4th order accuracy with the execution time equal to that of second order methods. Future research directions include application of this method to more realistic problems and development of a method that preserves both the total energy and total probability exactly.


\appendix

\section*{Appendix}

The target problem of the present paper stems from 
the quantum dynamics simulation of electronic wavefunction,
so as to understand the electrical conductivity of organic semiconductor materials. 
Organic semiconductor materials form the foundation of flexible devices \cite{hammock2013}.
The conductivity is simulated by the time evolution of  
a charged \lq carrier'
with a wavefunction $\Psi(\bm{r}, t)$,
for example, as in \cite{troisi2006}.
The carrier is classified into hole or electron,
where a hole or electron is positively or negatively charged, respectively. 
The wavepacket is expressed  
by the linear combination of $m$ given basis functions $\{ \chi_i(\bm{r}) \}_{i=1}^m$:
\begin{eqnarray}
\Psi(\bm{r}, t) = \sum_{i=1}^m u_i(t) \chi_i(\bm{r}),
\end{eqnarray}
with a complex coefficient vector of $\bm{u}(t) \equiv (u_1(t), u_2(t),\ldots, u_m(t))^{\top}$.
Hereafter we call $m$-site model,
where a site represents, typically, a molecule or atom.
The $i$-th component $u_i$ ($i=1,2, \ldots, m$) indicates the amplitude on the $i$-th site
and we assume that the weight of the carrier on the $i$-th site is
written as $n_i \equiv |u_i|^2$ ($\sum_{i=1}^m n_i =1$).
The simplest theoretical model is written as a linear problem of the form:
\begin{eqnarray}
 i \frac{d \bm{u}}{dt} =  K \bm{u}
\end{eqnarray}
with a given Hermitian matrix $K$. 
A matrix element $K_{ij}$ represents the transfer effect of the carrier 
between the $i$-th and $j$-th sites.
Since the transfer occurs only between the neighboring sites, the matrix $K$ is sparse. 
The energy 
\begin{eqnarray}
U_K \equiv \sum_{k,j} \bar{u}_k K_{kj} u_j 
    = \bm{u}^{\mathsf H} K \bm{u}. 
\end{eqnarray}
is conserved.
The time-evolution can be carried out by standard methods,
like the Crank-Nicolson method.
One of the authors (T. H.) carried out
quantum dynamics simulations using this linear model, see \cite{imachi2016,hoshi2016}, for example.

Recently, an advanced theoretical model was proposed \cite{terao2013,tada2018} by adding a non-linear energy term
\begin{eqnarray}
U_I \equiv U_I(n_1,n_2, \ldots, n_m) 
    =  - \frac{1}{2} \sum_{i=1}^m \gamma_i n_i^2.
\end{eqnarray}
to the linear model. In this model, the energy  $G \equiv U_K+U_I$ is conserved. 
The differential equation can be written as
\begin{eqnarray}
i \frac{d \bm{u}}{dt} =  (K+V_I) \bm{u},
\end{eqnarray}
where $V_I \equiv {\rm diag}(v_1, v_2, \ldots , v_m)$ and
$v_i \equiv - \gamma_i n_i$.
In addition, a $\delta$-function like term $V_E$ appears, when an impurity site is included in material
and the impurity can trap the carrier. 
A proper algorithm for the NLS-type equation with $\delta$-function like terms 
is not trivial and the application researchers would like to compare them. 
In general, a numerically robust algorithm can adopt a large time interval $h$ and saves the iteration steps $T/h$ for a given time period for the simulation $T$.
In particular, 
the reliable algorithm for the energy conservation is important, 
which motivates the present topic.

A future aspect for real researches is 
to construct an adaptive solver which contains
many integrator algorithms, like those in Figs.~\ref{fig:fig1} and \ref{fig:fig2}.
One will obtain automatically the optimal algorithm and the optimal
time interval $\delta t$
during the simulations,
when a proper evaluation for the time cost is realized.


\end{document}